\documentclass{iopconfser}

\usepackage{cite}

\usepackage{amsmath,amssymb,amsthm}
\usepackage{latexsym}
\usepackage{stmaryrd}
\usepackage{wasysym}
\usepackage{multicol}
\usepackage{multirow}
\usepackage{dirtytalk}
\usepackage{url}
\usepackage{fontawesome} 

\newtheorem{theor}{Theorem}
\newtheorem{claim}[theor]{Claim}

\theoremstyle{definition}
\newtheorem{simp}{Simplification}

\newtheorem{prop}[theor]{Proposition}
\newtheorem{lem}[theor]{Lemma}
\newtheorem{cor}[theor]{Corollary}
\newtheorem{deff}{Definition}

\newtheorem*{prob}{Problem formulation}
\newtheorem*{comment}{Comment}

\newtheorem{open}{Open problem}
\newtheorem{ex}{Example}
\newtheorem{nonex}[ex]{Non-example}
\theoremstyle{remark}
\newtheorem{rem}{Remark}

\def\BB{\mathbb}

\usepackage{bm}

\begin{document}

\pagestyle{plain} 

\title{Kontsevich graphs act on Nambu--Poisson brackets, II. \\ The tetrahedral flow is a coboundary in 4D}

\author{Mollie S Jagoe Brown, Floor Schipper and  Arthemy V Kiselev}

\affil{Bernoulli Institute for Mathematics, Computer Science and Artificial Intelligence, University of Groningen, P.O. Box 407, 9700 AK Groningen, The Netherlands}

\email{m.s.jagoe.brown@gmail.com, f.m.schipper@rug.nl, a.v.kiselev@rug.nl}

\begin{abstract}
Kontsevich constructed a map from suitable cocycles in the graph complex to infinitesimal deformations of Poisson bi-vector fields. Under the deformations, the bi-vector fields remain Poisson. We ask, are these deformations trivial, meaning, do they amount to a change of coordinates along a vector field? We examine this question for the tetrahedron, the smallest nontrivial suitable graph cocycle in the Kontsevich graph complex, and for the class of Nambu--Poisson brackets on $\mathbb{R}^d$.

Within Kontsevich's graph calculus, we use dimension-specific micro-graphs, in which each vertex represents an ingredient of the Nambu--Poisson bracket. For the tetrahedron, Kontsevich knew that the deformation is trivial for $d=2$ (1996). In 2020, Buring and the third author found that the deformation is trivial for $d=3$. Building on these discoveries, we now establish that the deformation is trivial for $d=4$.
\end{abstract}

\section{Introduction}

Take an arbitrary Poisson geometry: energy is transformed into motion by a Poisson bracket. We ask: is this model isolated, or part of a larger family? To examine this, we deform the Poisson bracket. If the deformation simply amounts to a change of coordinates of the model, then the system is isolated with respect to that deformation. If the deformation causes a nontrivial change of the Poisson bracket, then we say that the system is one in a larger family of systems. The incoming energy has stayed the same, but the outgoing motion has changed, see Chapter 13 on Deformation Quantization in~\cite{poisson}.

Kontsevich constructed a map from `good graphs' $\gamma$ in the Kontsevich graph complex, a differential graded Lie algebra of non-directed graphs, to bi-vector field flows $\dot{P}=Q^{\gamma}_d(P)$, which express the deformation of the Poisson bi-vector field $P$ on $\mathbb{R}^d$ by $\gamma$, see \cite{ascona}. The graphs $\gamma$ are cocycles, meaning that for the differential $\sf{d}=[\bullet\!\text{\textbf{--}}\!\bullet,\cdot]$ in the Kontsevich graph complex, we have $$\sf{d}(\gamma)=[\bullet\!\text{\textbf{--}}\!\bullet,\gamma]=0.$$ The associated bi-vector field flow $\dot{P}=Q_d^{\gamma}(P)$ is a Poisson cocycle for the differential $\partial_P=\llbracket P,\cdot\rrbracket$, meaning $$\partial_P\big(Q_d^\gamma(P)\big)=\llbracket P, Q^{\gamma}_d(P)\rrbracket =0.$$ We enquire if $Q^{\gamma}_d(P)$ is a coboundary, which would mean there exists some trivialising vector field $\vec{X}_d^{\gamma}(P)$ such that $Q^{\gamma}_d(P)=\llbracket P, \vec{X}_d^{\gamma}(P)\rrbracket$, which implies that the deformation $\dot{P}=Q^{\gamma}_d(P)$ is trivial; the coordinates change along the vector field $\vec{X}_d^{\gamma}(P)$. We specifically deform the class of Nambu--Poisson brackets by the simplest nontrivial graph in the Kontsevich graph complex, the tetrahedron $\gamma_3$. In this paper, we present the solution $\vec{X}^{\gamma_3}_{d=4}(P(\varrho,\boldsymbol{a}))$ found in dimension $d=4$ for Nambu--Poisson brackets $P(\varrho,\boldsymbol{a})$.

The authors recommend the following reading order of proceedings, by names of first authors: Kiselev~\cite{avk}, Jagoe Brown (this paper), and finally, Schipper~\cite{fs}.

This paper is structured in the following way. In section \ref{preliminaries}, we introduce the preliminaries necessary to approach the problem, then formulate it. In section \ref{vf's}, we examine the solution to the problem in dimensions two, three, and finally, dimension four, while presenting a series of simplifications which were crucial to obtain the new solution in dimension four. In section \ref{discussion}, we discuss the up-down behaviour of the problem from one dimension to another, and finally conclude. 

\section{Preliminaries}\label{preliminaries}

\subsection{Basic concept}

The theory behind this problem is due to Kontsevich, and is applicable to any class of Poisson bracket on an affine manifold, in any dimension. Recall that we can express any Poisson bracket in terms of a bi-vector field, in the following way: $$\{f,g\}=P(f,g).$$ Deforming a Poisson bi-vector field $P$ by a suitable graph cocycle $\gamma$ in the Kontsevich graph complex is expressed as $$\dot{P}=Q^{\gamma}(P),$$ where $Q^{\gamma}(P)$ is an infinitesimal symmetry built of as many copies of $P$ as there are vertices in $\gamma$, see~\cite{rb}.

The setting of the problem is $\BB{R}^d$, with Cartesian coordinates given by $\BB{R}^d\ni\boldsymbol{x}$ $=(x_1,x_2,...,x_d)$. We deform the class of Nambu--Poisson brackets by the tetrahedron $\gamma_3$.

\begin{deff}[Nambu--Poisson bracket]\label{npdef}
The generalised Nambu-determinant Poisson bracket in dimension $d$ for two smooth functions $f,g\in C^\infty(\BB{R}^d)$ is given as $$\{f,g\}_d(\boldsymbol{x})=\varrho(\boldsymbol{x})\cdot \det \begin{pmatrix}
    f_{x_1} & g_{x_1} & a^1_{x_1} & a^2_{x_1} & \hdots & a^{d-2}_{x_1} \\
    f_{x_2} & g_{x_2} & a^1_{x_2} & a^2_{x_2} & \hdots & a^{d-2}_{x_2} \\
    \vdots  & \vdots & \vdots & \vdots  & & \vdots  \\
    f_{x_d} & g_{x_d} & a^1_{x_d} & a^2_{x_d} & \hdots & a^{d-2}_{x_d}
\end{pmatrix}(\boldsymbol{x}),$$
where $a^1,...,a^{d-2}\in C^\infty(\BB{R}^d)$ are Casimirs, which Poisson-commute with any function. The function $\varrho$ is the inverse density, or the coefficient of a $d$-vector field.
\end{deff}

\begin{deff}[The tetrahedron $\gamma_3$]
The tetrahedron $\gamma_3$ is the smallest nontrivial suitable graph cocycle of the Kontsevich graph complex, a differential graded Lie algebra of non-directed graphs. By $\gamma_3$ being a cocycle, we mean that it satisfies $\sf{d}(\gamma_3)=0$, where $\sf{d}=[\bullet\!\text{\textbf{--}}\!\bullet,\cdot]$ is the differential in the Kontsevich graph complex. The graph $\gamma_3$ is constructed on 4 vertices and 6 edges. 
\end{deff}

\begin{deff}[The $\gamma_3$-flow]
    The $\gamma_3$-flow $Q^{\gamma_3}(P)$ is a bi-vector field built with four copies of $P$. It is an infinitesimal symmetry of the Jacobi identity for $P$, and corresponds to the deformation of $P$ by $\gamma_3$. It is obtained from $\gamma_3$ via the orientation morphism described in~\cite{rb}; the formula is given in~\cite{f16}.
\end{deff} 

\begin{prob}
To inspect whether the $\gamma_3$-flow is Poisson-trivial, we investigate if $Q^{\gamma_3}_d(P)$ is a 2-coboundary. This is equivalent to establishing the existence of a trivialising vector field $\vec{X}^{\gamma_3}_d(P)$ such that \begin{equation}\label{maineq}\dot{P}=Q^{\gamma_3}_d(P)=\llbracket P, \vec{X}^{\gamma_3}_d(P)\rrbracket,\end{equation} where $d$ is the dimension and $\llbracket\cdot,\cdot\rrbracket$ is the Schouten bracket\footnote{The Schouten bracket $\llbracket \cdot, \cdot\rrbracket$ is a unique extension of the commutator $[\cdot,\cdot]$ on the space of vector fields to the space of polyvector fields. By definition, the Schouten bracket coincides
with the Lie bracket when evaluated on 1-vectors. When evaluated on $p$-vector $X$, $q$-vector $Y$ and $r$-vector $Z$, the Schouten bracket satisfies the equations $\llbracket X, Y\rrbracket = -(-1)^{(p-1)(q-1)}\llbracket Y, X\rrbracket$ and $\llbracket X, Y \wedge  Z\rrbracket =
\llbracket X, Y\rrbracket  \wedge Z + (-1)^{q(p-1)}Y \wedge \llbracket X, Z\rrbracket$, see~\cite{tetra16}.}. We wish to solve equation (\ref{maineq}) in $d=4$. 
\end{prob}

We used software package $\textsf{gcaops}$\footnote{\url{https://github.com/rburing/gcaops}} (\textbf{G}raph \textbf{C}omplex \textbf{A}ction \textbf{O}n \textbf{P}oisson \textbf{S}tructures) for SageMath by Buring. With it, we input graph encodings, from which we obtained Formality graphs and then their formulas. We created a linear algebraic system and solved it for coefficients in the linear combination of graphs that encode $\vec{X}^{\gamma_3}_d(P)$. We used the High Performance Computing cluster at the University of Groningen, H\'abr\'ok. All code which gave results in this paper can be found as additional material attached to this paper and to \cite{fs}.

\subsection{Notation}

We solve equation (\ref{maineq}) on the level of formulas, using Kontsevich's graph calculus to write them. For this, we introduce the graph language created by Kontsevich, commonly used in deformation quantisation. Its main convenience is that formulas change with the dimension, but pictures of graphs do not change. We specifically use this graph language for graphs built of wedges $\smash{\xleftarrow{L}\!\!\bullet\!\!\xrightarrow{R}}$, which are Poisson bi-vector fields. The directed edges are derivations which act on the content of vertices. To write the graph encodings up to and including dimension four, we use the following convention: 
\begin{multicols}{2}
\begin{itemize}
    \item 0 represents the sink,
    \item 1, 2, 3 represent Levi--Civita symbols,
    \item 4, 5, 6 represent Casimirs $a^1$,
    \item 7, 8, 9 represent Casimirs $a^2$.
\end{itemize}
\end{multicols} \noindent We denote by $\phi$ the map of Formality graphs to their formulas obtained by Kontsevich's graph language. 

\begin{ex}
    Let us take the following graphs $\Gamma_1$ and $\Gamma_2$.

    \begin{multicols}{2}\begin{itemize}\item[]The encoding of $\Gamma_1$ is (0,1 ; 1,3 ; 1,2).\item[]The encoding of $\Gamma_2$ is (0,2 ; 1,3 ; 1,2).\end{itemize}\end{multicols} 
    \vspace{5pt} 
    \begin{multicols}{2}\begin{itemize} \item[] The graph $\Gamma_1$ is: \hspace{20pt} \raisebox{0pt}[6mm][4mm]{\unitlength=1.2mm
\special{em:linewidth 0.8pt}
\linethickness{0.5pt}
\begin{picture}(50,50)(5,5)
\put(-5,-7){
\begin{picture}(50.00,24.00)
\put(10.00,10.00){\circle*{1}} 
\put(17.00,17.0){\circle*{1}} 
\put(3.00,17.0){\circle*{1}} 
\put(10.00,10.00){\vector(0,-1){7.30}} 
\put(17.00,17.00){\vector(-1,0){14.00}} 
\put(3.00,17.00){\vector(1,-1){6.67}} 
\bezier{30}(3,17)(6.67,13.67)(9.67,10.33)
%
%
\put(17,17){\vector(-1,-1){6.67}} 
%
\bezier{52}(17.00,17.00)(16.33,23.33)(10.00,24.00) 
\bezier{52}(10.00,24.00)(3.67,23.33)(3.00,17.00) 
\put(16.8,18.2){\vector(0,-1){1}}
\put(10,10){\line(1,-2){2}} 
\bezier{52}(10,10)(20,5)(12,6)
\put(10,10){\vector(-3,1){0}} 
\put(9.5,12){{\tiny$1$}}
\put(18,17){{\tiny$2$}}
\put(0,17){{\tiny$3$}}
\put(9.5,1){{\tiny$0$}}
\put(7,6){{\tiny$i_1$}}
\put(16,7){{\tiny$i_2$}}
\put(15,13){{\tiny$j_1$}}
\put(9.5,18){{\tiny$j_2$}}
\put(3,13){{\tiny$k_1$}}
\put(9.5,22){{\tiny$k_2$}}
\end{picture}
}\end{picture}} \item[] The graph $\Gamma_2$ is: \hspace{20pt} \raisebox{0pt}[6mm][4mm]{\unitlength=1.2mm
\special{em:linewidth 0.8pt}
\linethickness{0.5pt}
\begin{picture}(50,50)(5,5)
\put(-5,-7){
\begin{picture}(50.00,24.00)
\put(10.00,10.00){\circle*{1}} 
\put(17.00,17.0){\circle*{1}} 
\put(3.00,17.0){\circle*{1}} 
\put(10.00,10.00){\vector(0,-1){7.30}} 
\put(17.00,17.00){\vector(-1,0){14.00}} 
\put(3.00,17.00){\vector(1,-1){6.67}} 
\bezier{30}(3,17)(6.67,13.67)(9.67,10.33)
%
%
\put(17,17){\vector(-1,-1){6.67}} 
%
\bezier{52}(17.00,17.00)(16.33,23.33)(10.00,24.00) 
\bezier{52}(10.00,24.00)(3.67,23.33)(3.00,17.00) 
\put(16.8,18.2){\vector(0,-1){1}}
\put(10,10){\line(1,0){10}}
\bezier{52}(20,10)(27,10)(18,17)
\put(18,17){\vector(-1,1){0}}
\put(9.5,12){{\tiny$1$}}
\put(18,17){{\tiny$2$}}
\put(0,17){{\tiny$3$}}
\put(9.5,1){{\tiny$0$}}
\put(7,6){{\tiny$i_1$}}
\put(16,7){{\tiny$i_2$}}
\put(15,13){{\tiny$j_1$}}
\put(9.5,18){{\tiny$j_2$}}
\put(3,13){{\tiny$k_1$}}
\put(9.5,22){{\tiny$k_2$}}
\end{picture}
}\end{picture}}

\end{itemize}\end{multicols} 

    \vspace{15pt}
    Let the dimension be two. The inert sums of the formulas for the graphs are as follows. $$\phi(\Gamma_1)=\sum_{\substack{i_1,i_2,\\j_1,j_2,\\k_1,k_2=1}}^{d=2} \varepsilon^{i_1i_2}\cdot \varepsilon^{j_1j_2}\cdot \varepsilon^{k_1k_2}\cdot \partial_{i_2j_1k_1}(\varrho)\cdot \partial_{k_2}(\varrho)\cdot \partial_{j_2}(\varrho)\cdot \partial_{i_1}(\text{ })$$
    $$\phi(\Gamma_2)=\sum_{\substack{i_1,i_2,\\j_1,j_2,\\k_1,k_2=1}}^{d=2} \varepsilon^{i_1i_2}\cdot \varepsilon^{j_1j_2}\cdot \varepsilon^{k_1k_2}\cdot \partial_{k_1j_1}(\varrho)\cdot \partial_{i_2k_2}(\varrho)\cdot \partial_{j_2}(\varrho)\cdot \partial_{i_1}(\text{ })$$
    The sums are constructed by taking the product of the content of vertices, which contain $\varrho$. The arrows act on vertices as derivations\footnote{We denote $\partial_i$ to mean the partial derivative with respect to $x_i$, represented by the arrow $i$; $\partial_{ij}$ is the partial derivative with respect to $x_i$ and $x_j$, so $\partial_{ij}=\partial_i\partial_j$.}. The Levi--Civita symbol encodes the determinant in the Nambu--Poisson bracket, see Definition \ref{npdef}. 
\end{ex}

\begin{deff}[The sunflower graph]
    A linear combination of the above Kontsevich graphs (graphs built of wedges $\smash{\xleftarrow{L}\!\!\bullet\!\!\xrightarrow{R}}$, see~\cite{avk,fs,rb}) can be expressed as the sunflower graph $$\text{sunflower }=\text{ }\raisebox{0pt}[6mm][4mm]{\unitlength=0.4mm
\special{em:linewidth 0.4pt}
\linethickness{0.4pt}
\begin{picture}(17,24)(5,5)
\put(-5,-7){
\begin{picture}(17.00,24.00)
\put(10.00,10.00){\circle*{1}}
\put(17.00,17.0){\circle*{1}}
\put(3.00,17.0){\circle*{1}}
\put(10.00,10.00){\vector(0,-1){7.30}}
\put(17.00,17.00){\vector(-1,0){14.00}}
\put(3.00,17.00){\vector(1,-1){6.67}}
\bezier{30}(3,17)(6.67,13.67)(9.67,10.33)
%
%
\put(17,17){\vector(-1,-1){6.67}}
\bezier{30}(17,17)(13.67,13.67)(10.33,10.33)
\bezier{52}(17.00,17.00)(16.33,23.33)(10.00,24.00)
\bezier{52}(10.00,24.00)(3.67,23.33)(3.00,17.00)
\put(16.8,18.2){\vector(0,-1){1}}
\put(10,17){\oval(18,18)}
\put(10,10){\line(1,0){10}}
\bezier{52}(20,10)(27,10)(21,16)
\put(21,16){\vector(-1,1){0}}
\end{picture}
}\end{picture}}=1\cdot\Gamma_1+2\cdot\Gamma_2.$$ The outer circle means that the outgoing arrow acts on the three vertices via the Leibniz rule. When the arrow acts on the upper two vertices, we obtain two isomorphic graphs, hence the coefficient 2 in the linear combination.
\end{deff}

\section{Vector fields trivialising the $\gamma_3$-flow in 2D, 3D, and now, in 4D}\label{vf's}

\subsection{The trivialising vector field for $\gamma_3$-flow in 2D expressed by Kontsevich graphs}

In 2D, any Poisson bracket is a Nambu--Poisson bracket because it is given by $$\{f,g\}_{d=2}(x,y)=\varrho(x,y)\cdot\det\begin{pmatrix}
    f_x&g_x\\f_y&g_y
\end{pmatrix}(x,y), \text{ that is,} \quad \{f,g\}_{d=2}(x,y)=\sum_{i,j=1}^{d=2}\varepsilon^{ij}\cdot\varrho\cdot\partial_i(f)\cdot\partial_j(g),$$ for some $\varrho$, where $i,j$ are indices, $\varepsilon^{ij}$ is the Levi--Civita symbol, and the Cartesian coordinates are expressed as $x_1=x$, $x_2=y$.

\begin{prop}[Cf. \cite{ascona}, \cite{tetra16}]\label{2D}
   There exists a trivialising vector field in 2D for the $\gamma_3$-flow. It is given by the sunflower graph $$\vec{X}^{\gamma_3}_{d=2}(P)= \phi(\text{sunflower}).$$ The sunflower gives a formula in 2D to solve equation (\ref{maineq}), namely $\dot{P}=Q^{\gamma_3}_{d=2}(P)=\llbracket P,\vec{X}_{d=2}^{\gamma_3}(P)\rrbracket.$ The affine space of solutions on graphs is of dimension 1; that is, the trivialising vector field $\vec{X}^{\gamma_3}_{d=2}(P)$ is unique up to 1-dimensional shifts, themselves encoded by Kontsevich graphs. A full descriptions of these shifts using graphs can be found in~\cite{fs}. 
\end{prop}

\begin{rem}
    This solution, the sunflower, when extended to higher dimensions $d>2$, is not a solution.
\end{rem}

\subsection{The trivialising vector field for $\gamma_3$-flow in 3D expressed by Nambu micro-graphs}

In 3D, with Cartesian coordinates expressed as $x_1=x$, $x_2=y$, $x_3=z$, the Nambu--Poisson bracket is given by $$\{f,g\}_{d=3}(x,y,z)=\varrho(x,y,z)\cdot\det\begin{pmatrix}
    f_x&g_x&a^1_x\\f_y&g_y&a^1_y\\f_z&g_z&a^1_z
\end{pmatrix}(x,y,z),$$ that is, $$\{f,g\}_{d=3}(x,y,z)=\sum_{i,j,k=1}^{d=3}\varepsilon^{ijk}\cdot \varrho\cdot \partial_i(f)\cdot\partial_j(g)\cdot\partial_k(a^1),$$ for some $\varrho$, where $i,j,k$ are indices and $\varepsilon^{ijk}$ is the Levi--Civita symbol which encodes the determinant.

\begin{deff}[Nambu micro-graph]
    Nambu micro-graphs are built using the Nambu--Poisson bracket $P(\varrho,\boldsymbol{a})$ as subgraphs with ordered and directed edges. The vertex of the source of each $P(\varrho,\boldsymbol{a})$ in dimension $d$ contains $\varepsilon^{i_1\hdots i_d}\varrho$, with $d$-many outgoing edges. The first two edges act on the arguments of that bi-vector field subgraph, and the last $d-2$ edges go to the Casimirs $a^1,...,a^{d-2}$.
\end{deff}

\begin{ex}
    Replacing the wedges of the sunflower graph with Nambu building blocks of any dimension gives an example of a Nambu micro-graph. 
\end{ex}

\begin{nonex}
    The encoding (0,1 ; 1,3,5 ; 2) is not a Nambu micro-graph. This can immediately be seen by the fact that the vertices do not have the same number of outgoing arrows. 
\end{nonex}

Instead of searching for a trivialising vector field $\vec{X}^{\gamma_3}_{d=3}(\varrho,\boldsymbol{a})$ over all 1-vector micro-graphs built of 3 Levi--Civita symbols and 3 Casimirs $a^1$, we restrict our scope. 

\begin{simp}
 Search for solutions $\vec{X}^{\gamma_3}_{d=3}(P(\varrho,\boldsymbol{a}))$ over the formulas given by Nambu micro-graphs.
\end{simp} 

\begin{deff}[$d$-descendants] The $d$-descendants of a given $(d'=2)$-dimensional Kontsevich graph is the set of Nambu micro-graphs obtained in the following way. Take a $(d'=2)$-dimensional Kontsevich graph. To each vertex, add $(d-2)$-many Casimirs by $(d-2)$-many outgoing edges. Extend the original incoming arrows to work via the Leibniz rule over the newly added Casimirs \footnote{The case of $d'>2$ is more delicate. This example gives all 4D-descendants from a given 3D Nambu micro-graph. We take the 3D Nambu micro-graph given by (0,1,4 ; 1,3,5 ; 1,2,6). Its 4D-descendants are given by $\sum_{i_1,i_2\in\{1,4,7\},j\in\{3,6,9\},k\in\{2,5,8\}}$(0,1,4,7 ; $i_1$,$j$,5,8 ; $i_2$,$k$,6,9). It is left to the reader to produce the generalised definition of $d$-descendants from $d'$-dimension.}.
\end{deff}

\begin{ex}
    We give the encodings of the 3D-descendants of the 2D sunflower. Recall that the sink is denoted by 0, and the Levi--Civita symbols by 1, 2, 3. In 3D, we have the three Casimirs $a^1$ denoted by 4, 5, 6. Then, $$\text{sunflower }=\text{ }\raisebox{0pt}[6mm][4mm]{\unitlength=0.4mm
\special{em:linewidth 0.4pt}
\linethickness{0.4pt}
\begin{picture}(17,24)(5,5)
\put(-5,-7){
\begin{picture}(17.00,24.00)
\put(10.00,10.00){\circle*{1}}
\put(17.00,17.0){\circle*{1}}
\put(3.00,17.0){\circle*{1}}
\put(10.00,10.00){\vector(0,-1){7.30}}
\put(17.00,17.00){\vector(-1,0){14.00}}
\put(3.00,17.00){\vector(1,-1){6.67}}
\bezier{30}(3,17)(6.67,13.67)(9.67,10.33)
%
%
\put(17,17){\vector(-1,-1){6.67}}
\bezier{30}(17,17)(13.67,13.67)(10.33,10.33)
\bezier{52}(17.00,17.00)(16.33,23.33)(10.00,24.00)
\bezier{52}(10.00,24.00)(3.67,23.33)(3.00,17.00)
\put(16.8,18.2){\vector(0,-1){1}}
\put(10,17){\oval(18,18)}
\put(10,10){\line(1,0){10}}
\bezier{52}(20,10)(27,10)(21,16)
\put(21,16){\vector(-1,1){0}}
\end{picture}
}\end{picture}}= 1\cdot\Gamma_1+2\cdot\Gamma_2=1\cdot \text{(0,1 ; 1,3 ; 1,2)}+2\cdot\text{(0,2 ; 1,3 ; 1,2)}$$ gives $$\text{3D-descendants }=\sum_{\substack{i_1,i_2\in\{1,4\}\\j\in\{3,6\}\\k\in\{2,5\}}} \text{(0,1,4 ; $i_1$,$j$,5 ; $i_2$,$k$,6) }+\sum_{\substack{i_1,i_2\in\{1,4\}\\j\in\{3,6\}\\k_1,k_2\in\{2,5\}}}\text{(0,$k_1$,4 ; $i_1$,$j$,5 ; $i_2$,$k_2$,6)}.$$
\end{ex}

\begin{lem}
    There are 41 distinct micro-graphs in the 3D expansion of the 2D sunflower, their encodings are below. The first 10 come from $\Gamma_1$, the next 31 come from $\Gamma_2$. In bold are those whose Formality graphs $\Gamma$ give formulas equal to zero $\phi(\Gamma)=0$. A graph is \textit{zero} when it has a symmetry under which it is skew. 
    \small
\begin{multicols}{3}
\begin{itemize}
\item[1.](0,1,4 ; 1,3,5 ; 1,2,6)    
\item[2.](0,1,4 ; 4,3,5 ; 4,2,6)
\item[3.](0,1,4 ; 1,6,5 ; 1,5,6)
\item[4.](0,1,4 ; 4,6,5 ; 4,5,6)
\item[5.](0,1,4 ; 1,3,5 ; 4,2,6)
\item[6.](0,1,4 ; 1,6,5 ; 1,2,6)
\item[7.](0,1,4 ; 1,6,5 ; 4,2,6)
\item[8.](0,1,4 ; 4,6,5 ; 1,2,6)
\item[9.](0,1,4 ; 4,6,5 ; 4,2,6)
\item[10.]\textbf{(0,1,4 ; 1,6,5 ; 4,5,6)}\hspace{0.3cm}$\uparrow\Gamma_1$

\vspace{0.05cm}
      \hrule
      \vspace{-0.05cm}
\item[11.](0,2,4 ; 1,3,5 ; 1,2,6)\hspace{0.73cm}$\downarrow\Gamma_2$
\item[12.](0,2,4 ; 4,3,5 ; 1,2,6)
\item[13.]\textbf{(0,2,4 ; 1,3,5 ; 4,2,6)}
\item[14.](0,2,4 ; 4,3,5 ; 4,2,6)
\item[15.](0,2,4 ; 1,6,5 ; 1,2,6)
\item[16.](0,2,4 ; 4,6,5 ; 1,2,6)
\item[17.](0,2,4 ; 1,6,5 ; 4,2,6)
\item[18.](0,2,4 ; 4,6,5 ; 4,2,6)
\item[19.](0,5,4 ; 1,3,5 ; 1,2,6)
\item[20.]\textbf{(0,5,4 ; 4,3,5 ; 1,2,6)}
\item[21.]\textbf{(0,5,4 ; 1,3,5 ; 4,2,6)} 
\item[22.](0,5,4 ; 4,3,5 ; 4,2,6)
\item[23.](0,5,4 ; 1,6,5 ; 1,2,6)
\item[24.]\textbf{(0,5,4 ; 4,6,5 ; 1,2,6)}
\item[25.]\textbf{(0,5,4 ; 1,6,5 ; 4,2,6)}
\item[26.](0,5,4 ; 4,6,5 ; 4,2,6)
\item[27.](0,2,4 ; 1,3,5 ; 1,5,6)
\item[28.](0,2,4 ; 4,3,5 ; 1,5,6)
\item[29.]\textbf{(0,2,4 ; 1,3,5 ; 4,5,6)}
\item[30.](0,2,4 ; 4,3,5 ; 4,5,6)
\item[31.](0,2,4 ; 1,6,5 ; 1,5,6)
\item[32.]\textbf{(0,2,4 ; 4,6,5 ; 1,5,6)}
\item[33.]\textbf{(0,2,4 ; 1,6,5 ; 4,5,6)}
\item[34.](0,2,4 ; 4,6,5 ; 4,5,6)
\item[35.](0,5,4 ; 1,3,5 ; 1,5,6)
\item[36.](0,5,4 ; 4,3,5 ; 1,5,6)
\item[37.]\textbf{(0,5,4 ; 1,3,5 ; 4,5,6)}
\item[38.]\textbf{(0,5,4 ; 4,3,5 ; 4,5,6)} \textit{zero}
\item[39.](0,5,4 ; 1,6,5 ; 1,5,6)
\item[40.](0,5,4 ; 4,6,5 ; 1,5,6)
\item[41.]\textbf{(0,5,4 ; 4,6,5 ; 4,5,6)} \textit{zero}
\end{itemize}
\end{multicols} 
\end{lem}

Our next simplification is a lucky guess, in contrast with simplification 1.

\begin{simp} Search for the trivialising vector field $\vec{X}^{\gamma_3}_{d=3}(P)$ over 41 3D-descendants of the 2D sunflower, from the above Lemma. 
\end{simp}

\begin{cor}
    Simplifications 1 and 2 make the problem smaller: $$366\xrightarrow{\#1,\#2}41,$$ while still allowing us to reach a solution. Here, 366 is the number of all 1-vector micro-graphs built of 3 Levi--Civita symbols and 3 Casimirs $a^1$; 41 is the number of 3D-descendants of the 2D sunflower.
\end{cor}


\begin{prop}\label{3D} There exists a trivialising vector field $\vec{X}^{\gamma_3}_{d=3}(P)=\phi(X^{\gamma_3}_{d=3})$ in 3D. It is given over 10 3D-descendants of the 2D sunflower:\begin{align*}
    X^{\gamma_3}_{d=3} &= 8\cdot(0,1,4 ; 1,3,5 ; 1,2,6) +24\cdot (0,1,4 ; 1,6,5 ; 4,2,6)+8\cdot (0,1,4 ; 4,3,5 ; 4,2,6)\\ &+24\cdot (0,1,4 ; 4,6,5 ; 4,2,6)+12\cdot (0,1,4 ; 4,6,5 ; 4,5,6)+16\cdot(0,2,4 ; 1,3,5 ; 1,2,6)\\&+16\cdot (0,2,4 ; 1,3,5 ; 1,5,6)+12\cdot(0,2,4 ; 1,6,5 ; 1,5,6)+16\cdot(0,2,4 ; 4,3,5 ; 1,5,6)+24\cdot(0,5,4 ; 1,3,5 ; 1,2,6).
\end{align*}This is a linear combination of Nambu micro-graphs which gives a formula in 3D to solve equation (\ref{maineq}), namely $\dot{P}=Q^{\gamma_3}_{d=3}(P)=\llbracket P,\vec{X}^{\gamma_3}_{d=3}(P)\rrbracket.$ The affine space of solutions on graphs is of dimension 3; that is, the trivialising vector field $\vec{X}^{\gamma_3}_{d=3}(P)$ is unique up to a 3-dimensional space of Poisson 1-cocycles $X$ with $\llbracket P, X\rrbracket = 0$, where $X$ is encoded by Nambu micro-graphs. A full description of these shifts using graphs can be found in~\cite{fs}.
\end{prop} 

We verified that the deformation of $P$ by $\gamma_3$ is trivial in 3D: it amounts to a change of coordinates. The space of 3D-descendants from the 2D sunflower is sufficient to find a solution in 3D. We now ask: \begin{quote}\textbf{Can we find a 4D solution over the space of 4D-descendants from the 3D trivialising vector field? (Answer: no! But from the 2D sunflower: yes!)}\end{quote}

\subsection{New result: the trivialising vector field for $\gamma_3$-flow in 4D}

In 4D, with Cartesian coordinates expressed as $x_1=x$, $x_2=y$, $x_3=z$, $x_4=w$, the Nambu--Poisson bracket is given by $$\{f,g\}_{d=4}(x,y,z,w)=\varrho(x,y,z,w)\cdot\det\begin{pmatrix}
    f_x&g_x&a^1_x&a^2_x\\f_y&g_y&a^1_y&a^2_y\\f_z&g_z&a^1_z&a^2_z\\f_w&g_w&a^1_w&a^2_w
\end{pmatrix}(x,y,z,w),$$ that is, $$\{f,g\}_{d=4}(x,y,z,w)=\sum_{i,j,k,\ell=1}^{d=4}\varepsilon^{ijk\ell}\cdot \varrho\cdot \partial_i(f)\cdot\partial_j(g)\cdot\partial_k(a^1)\cdot\partial_{\ell}(a^2),$$ for some $\varrho$, where $i,j,k,\ell$ are indices and $\varepsilon^{ijk\ell}$ is the Levi--Civita symbol which encodes the determinant.

To tackle the full problem was non-viable, see Appendix 1 in~\cite{skew23} for details. The main issue is that the size of the problem increases with the dimension: $\phi(Q^{\gamma_3}_d(P))$ is 1 line in 2D, 2 pages in 3D, 3GB in~4D.

We begin to apply the simplifications that we have introduced so far. 
\begin{itemize} 
\item Simplification 1: only look over graphs built with the 4D Nambu building blocks $P(\varrho,a^1,a^2)$. The vertex of the source of the building blocks contains $\varepsilon^{ijk\ell}\varrho$.

\item Simplification 2: only look over 4D-descendants of the 2D sunflower. These are 324 graphs, which give 123 linearly independent formulas; their encodings can be found in the Appendix.

\item[\faWarning] There are two Casimirs $a^1,a^2$, yielding an extra property to take into account.
\end{itemize}

\begin{lem} The Nambu--Poisson bracket $P(\varrho,a^1,a^2)$ is skew-symmetric under the swap $a^1\rightleftarrows a^2$: $$P(\varrho,a^1,a^2)=-P(\varrho,a^2,a^1).$$ The $\gamma_3$-flow $Q^{\gamma_3}_{d=4}(P)$ is built of four copies of $P$, therefore $Q^{\gamma_3}_{d=4}(P)$ is symmetric under $a^1\rightleftarrows a^2$; by swapping $a^1$ and $a^2$, we accumulate four minus signs: \begin{multline*}Q^{\gamma_3}_{d=4}\Big(P(\varrho,a^2,a^1)\otimes P(\varrho,a^2,a^1)\otimes P(\varrho,a^2,a^1)\otimes P(\varrho,a^2,a^1)\Big)\\=(-)^4Q^{\gamma_3}_{d=4}\Big(P(\varrho,a^1,a^2)\otimes P(\varrho,a^1,a^2)\otimes P(\varrho,a^1,a^2)\otimes P(\varrho,a^1,a^2)\Big)\\=Q^{\gamma_3}_{d=4}\Big(P(\varrho,a^1,a^2)\otimes P(\varrho,a^1,a^2)\otimes P(\varrho,a^1,a^2)\otimes P(\varrho,a^1,a^2)\Big).\end{multline*} Therefore, to find a vector field $\vec{X}^{\gamma_3}_{d=4}(P)$ such that $$\dot{P}=Q^{\gamma_3}_{d=4}(P)=\llbracket P,\vec{X}^{\gamma_3}_{d=4}(P)\rrbracket,$$ we need to find $\vec{X}^{\gamma_3}_{d=4}(P)$ which is skew-symmetric under $a^1\rightleftarrows a^2$. This can be seen by the fact that $\vec{X}^{\gamma_3}_{d=4}(P)$ is built of three copies of $P$, so accumulates three minus signs when swapping $a^1$ and $a^2$, therefore gives $(-)^3=-$, therefore is skew-symmetric under $a^1\rightleftarrows a^2$.
\end{lem}

To take this into account, we use the 324 4D-descendants of the 2D sunflower obtained by simplifications 1 and 2, identify the 123 ones with linearly independent formulas, and skew-symmetrise them. That is, for each $4D$-descendant $\Gamma$ of the 2D sunflower we construct a skew pair: $$\text{skew pair}=\tfrac{1}{2}\Bigl(\phi\bigl(\Gamma(a^1,a^2)\bigr)-\phi\bigl(\Gamma(a^2,a^1)\bigr)\Bigr).$$ To construct skew pairs, we take the formula of the graph $\Gamma$ with ordering of edges to Casimirs $a^1,a^2$ with $a^1\prec a^2$, and subtract the formula of the graph $\Gamma$ with ordering $a^2\prec a^1$. We divide by 2 to preserve the coefficients. By construction, each skew pair is purely obtained at the level of formulas\footnote{See \url{https://github.com/rburing}.}.

\begin{lem}
    There are 64 linearly independent skew pairs, see Appendix.
\end{lem}

\begin{simp}
    Search over these 64 skew pairs for a trivialising vector field $\vec{X}^{\gamma_3}_{d=4}(P)$.
\end{simp}

\begin{cor}\label{3simps}
    The three simplifications 1, 2, and 3 reduced the size of our problem 300 times: $$19\text{ }957 \xrightarrow{\#1,\#2} 324 \xrightarrow{\#3} 64.$$ Here, 19 957 is the number of all 1-vector micro-graphs built of 4 Levi--Civita symbols, 4 Casimirs $a^1$ and 4 Casimirs $a^2$; 324 is the number of 4D-descendants of the 2D sunflower; 64 is the number of skew pairs obtained from the 123 linearly independent formulas of the 324 4D-descendants. 
\end{cor}

\begin{prop}\label{4D} There exists a trivialising vector field $\vec{X}_{d=4}^{\gamma_3}(P)=\phi(X_{d=4}^{\gamma_3})$ in 4D. Searching over the 64 skew pairs on the High Performing Computing cluster H\'abr\'ok took 10 hours. It is given over 27 skew pairs of 1-vector Nambu micro-graphs:
\tiny \begin{multicols}{2}\begin{itemize}
\item[${X}^{\gamma_3}_{d=4}$]= $-8\cdot\Big($(0,1,4,7 ; 1,3,5,8 ; 1,2,6,9) - (0,1,7,4 ; 1,3,8,5 ; 1,2,9,6)$\Big)$
\item[+]$-48\cdot\Big($(0,1,4,7 ; 1,6,5,8 ; 4,2,6,9) - (0,1,7,4 ; 1,9,8,5 ; 7,2,9,6)$\Big)$
\item[+]$-16\cdot\Big($(0,1,4,7 ; 4,3,5,8 ; 4,2,6,9) - (0,1,7,4 ; 7,3,8,5 ; 7,2,9,6)$\Big)$
\item[+]$-48\cdot\Big($(0,1,4,7 ; 4,6,5,8 ; 4,2,6,9) - (0,1,7,4 ; 7,9,8,5 ; 7,2,9,6)$\Big)$
\item[+]$-48\cdot\Big($(0,1,4,7 ; 4,9,5,8 ; 4,2,6,9) - (0,1,7,4 ; 7,6,8,5 ; 7,2,9,6)$\Big)$
\item[+]$-16\cdot\Big($(0,1,4,7 ; 4,3,5,8 ; 7,2,6,9) - (0,1,7,4 ; 7,3,8,5 ; 4,2,9,6)$\Big)$
\item[+]$-48\cdot\Big($(0,1,4,7 ; 4,6,5,8 ; 7,2,6,9) - (0,1,7,4 ; 7,9,8,5 ; 4,2,9,6)$\Big)$
\item[+]$12\cdot\Big($(0,1,4,7 ; 1,6,5,8 ; 1,5,6,9) - (0,1,7,4 ; 1,9,8,5 ; 1,8,9,6)$\Big)$
\item[+]$-24\cdot\Big($(0,1,4,7 ; 4,6,5,8 ; 4,5,6,9) - (0,1,7,4 ; 7,9,8,5 ; 7,8,9,6)$\Big)$
\item[+]$24\cdot\Big($(0,1,4,7 ; 4,9,5,8 ; 7,5,6,9) - (0,1,7,4 ; 7,6,8,5 ; 4,8,9,6)$\Big)$
\item[+]$-24\cdot\Big($(0,1,4,7 ; 7,6,5,8 ; 7,5,6,9) - (0,1,7,4 ; 4,9,8,5 ; 4,8,9,6)$\Big)$
\item[+]$-16\cdot\Big($(0,2,4,7 ; 1,3,5,8 ; 1,2,6,9) - (0,2,7,4 ; 1,3,8,5 ; 1,2,9,6)$\Big)$
\item[+]$-32\cdot\Big($(0,2,4,7 ; 1,3,5,8 ; 1,5,6,9) - (0,2,7,4 ; 1,3,8,5 ; 1,8,9,6)$\Big)$
\item[+]$-12\cdot\Big($(0,2,4,7 ; 1,6,5,8 ; 1,5,6,9) - (0,2,7,4 ; 1,9,8,5 ; 1,8,9,6)$\Big)$
\item[+]$48\cdot\Big($(0,2,4,7 ; 1,9,5,8 ; 4,5,6,9) - (0,2,7,4 ; 1,6,8,5 ; 7,8,9,6)$\Big)$
\item[+]$-48\cdot\Big($(0,2,4,7 ; 1,9,5,8 ; 7,5,6,9) - (0,2,7,4 ; 1,6,8,5 ; 4,8,9,6)$\Big)$
\item[+]$-32\cdot\Big($(0,2,4,7 ; 4,3,5,8 ; 1,5,6,9) - (0,2,7,4 ; 7,3,8,5 ; 1,8,9,6)$\Big)$
\item[+]$-32\cdot\Big($(0,2,4,7 ; 7,3,5,8 ; 1,5,6,9) - (0,2,7,4 ; 4,3,8,5 ; 1,8,9,6)$\Big)$
\item[+]$-48\cdot\Big($(0,5,4,7 ; 1,3,5,8 ; 1,2,6,9) - (0,8,7,4 ; 1,3,8,5 ; 1,2,9,6)$\Big)$
\item[+]$-48\cdot\Big($(0,5,4,7 ; 1,9,5,8 ; 1,2,6,9) - (0,8,7,4 ; 1,6,8,5 ; 1,2,9,6)$\Big)$
\item[+]$-96\cdot\Big($(0,5,4,7 ; 1,9,5,8 ; 4,2,6,9) - (0,8,7,4 ; 1,6,8,5 ; 7,2,9,6)$\Big)$
\item[+]$-48\cdot\Big($(0,5,4,7 ; 4,9,5,8 ; 4,2,6,9) - (0,8,7,4 ; 7,6,8,5 ; 7,2,9,6)$\Big)$
\item[+]$-48\cdot\Big($(0,5,4,7 ; 7,9,5,8 ; 4,2,6,9) - (0,8,7,4 ; 4,6,8,5 ; 7,2,9,6)$\Big)$
\item[+]$-48\cdot\Big($(0,5,4,7 ; 7,9,5,8 ; 7,2,6,9) - (0,8,7,4 ; 4,6,8,5 ; 4,2,9,6)$\Big)$
\item[+]$24\cdot\Big($(0,5,4,7 ; 1,6,5,8 ; 1,5,6,9) - (0,8,7,4 ; 1,9,8,5 ; 1,8,9,6)$\Big)$
\item[+]$48\cdot\Big($(0,5,4,7 ; 4,6,5,8 ; 7,5,6,9) - (0,8,7,4 ; 7,9,8,5 ; 4,8,9,6)$\Big)$
\item[+]$48\cdot\Big($(0,5,4,7 ; 4,9,5,8 ; 7,5,6,9) - (0,8,7,4 ; 7,6,8,5 ; 4,8,9,6)$\Big)$.
    \end{itemize}
    \end{multicols}
\normalsize
This is a linear combination of skew pairs which gives a formula in 4D to solve equation (\ref{maineq}), namely $\dot{P}=Q^{\gamma_3}_{d=4}(P)=\llbracket P,\vec{X}^{\gamma_3}_{d=4}(P)\rrbracket.$ The affine space of solutions on graphs is of dimension 7; that is, the trivialising vector field $\vec{X}^{\gamma_3}_{d=4}(P)$ is unique up to a 7-dimensional space of Poisson 1-cocycles $X$ with $\llbracket P, X\rrbracket = 0$, where $X$ is encoded by Nambu micro-graphs. A full description of these shifts using graphs can be found in~\cite{fs}.
\end{prop}

As far as we understand, nothing could have predicted this result. Without this series of simplifications (see Corollary \ref{3simps}), approaching the problem in dimension 4 was impossible in~\cite{skew23} two years ago.

\section{Discussion}\label{discussion}

We summarise in which sense the graphs used in problems 2D, 3D, and 4D were different. In 2D, we found a solution over Kontsevich graphs; in 3D over Nambu micro-graphs; in 4D over skew pairs, obtained by skew-symmetrising formulas of Nambu micro-graphs. The problem has been solved at the level of formulas; the graphs provide a roadmap to find the few ones which have given solutions so far. Each dimension has presented a peculiarity related to properties of the class of Nambu--Poisson brackets. In 3D, we encoded the first instance where the Nambu--Poisson bracket is degenerate (of rank $\leq2$), that is, different than the maximal-rank Poisson bracket. In 4D, we skew-symmetrised the formulas of micro-graphs to account for solutions' anti-symmetry with respect to two Casimirs. Overall, we notice an interplay between expansions of graphs between dimensions, and the `goodness' of their formulas. The behaviour of solutions in dimension $d$ to dimensions $d-1$ and $d+1$ exhibits two curious properties.
\begin{claim}
    We can project a 4D solution down to a 3D solution by setting the last Casimir equal to the last coordinate: $a^2=w$. Similarly, we can project a 3D solution down to a 2D solution by setting the Casimir equal to the last coordinate: $a^1=z$. Formulas project down to previously found formulas: $$\phi\big(\vec{X}^{\gamma_3}_{d=4}(P)\big)\xrightarrow{a^2=w} \phi\big(\vec{X}^{\gamma_3}_{d=3}(P)\big)\xrightarrow{a^1=z}\phi\big(\vec{X}^{\gamma_3}_{d=2}(P)\big).$$ 
\end{claim}
\begin{claim}
    There exist solutions in 3D and 4D over the descendants of a 2D solution, the sunflower. But descendants of known solutions in 3D do not give solutions in 4D.
\end{claim}
\begin{comment}
This would have been practical for reducing computing time. We have 41 3D graphs obtained from expanding the sunflower, and solutions in 3D over 10 such graphs. Moving up to a higher dimension, it would be ideal to search over descendants of a 3D solution, but we observe this is not possible.
\end{comment}

\begin{open}
    What is the formula of the trivialising vector field $\vec{X}^{\gamma_3}_d(P)$ in dimensions $d\geq5$, if it exists at all?
\end{open}

\begin{open}
    In which dimensions are deformations of Nambu--Poisson brackets by other `good graphs' in the Kontsevich graph complex, such as $\gamma_5,\gamma_7,[\gamma_3,\gamma_5]$, (non)trivial? 
\end{open}

The order of the problem also increases with the choice of `good graph', for instance computing with $\gamma_5$ in dimension 3 is more costly than computing with $\gamma_3$ in dimension 4, see~\cite{skew21,avk}.

\section{Conclusion}

We observed that in dimension 4, the Kontsevich $\gamma_3$-flow is trivial for the class of Nambu--Poisson brackets. In other words, the deformation of the class of Nambu--Poisson brackets by $\gamma_3$ in 4D amounts to a change of coordinates along a vector field $\vec{X}^{\gamma_3}_{d=4}(P)$. So far, in dimensions 2, 3 and 4, we have that the Nambu--Poisson system is isolated in that it withstands the deforming action of the Kontsevich graph $\gamma_3$. Achieving this result was only possible by using the series of simplifications introduced here. 

\section*{Acknowledgements}

The authors thank the organizers of the ISQS28 conference on
1--5 July 2024 in CVUT Prague for a dynamic and welcoming atmosphere; and the Center for Information Technology of the University of Groningen for access to the High Performance Computing cluster, H\'abr\'ok. The authors are especially grateful to R.~Buring for computational support for his software package $\textsf{gcaops}$. The participation of M.S.~Jagoe Brown and F.~Schipper in the ISQS28 was
supported by the Master's Research Project funds at the Bernoulli Institute, University of Groningen; that of A.V.~Kiselev was supported by project~135110.


\newpage

\section*{Appendix}

There are 324 micro-graphs in the 4D expansion of the 2D sunflower, their encodings are given below. The first 81 come from $\Gamma_1$, the next 243 come from $\Gamma_2$. In bold are the 54 encodings whose Formality graphs $\Gamma$ give formulas equal to zero $\phi(\Gamma)=0$. A graph is zero when it has a symmetry under which it is skew. 
\tiny 
\begin{multicols}{3}
\begin{itemize}
    \item[1.] (0, 1, 4, 7, 1, 3, 5, 8, 1, 2, 6, 9)
\item[2.] (0, 1, 4, 7, 1, 6, 5, 8, 1, 2, 6, 9)
\item[3.] (0, 1, 4, 7, 1, 9, 5, 8, 1, 2, 6, 9)
\item[4.] (0, 1, 4, 7, 1, 3, 5, 8, 4, 2, 6, 9)
\item[5.] (0, 1, 4, 7, 1, 6, 5, 8, 4, 2, 6, 9)
\item[6.] (0, 1, 4, 7, 1, 9, 5, 8, 4, 2, 6, 9)
\item[7.] (0, 1, 4, 7, 1, 3, 5, 8, 7, 2, 6, 9)
\item[8.] (0, 1, 4, 7, 1, 6, 5, 8, 7, 2, 6, 9)
\item[9.] (0, 1, 4, 7, 1, 9, 5, 8, 7, 2, 6, 9)
\item[10.] (0, 1, 4, 7, 4, 3, 5, 8, 1, 2, 6, 9)
\item[11.] (0, 1, 4, 7, 4, 6, 5, 8, 1, 2, 6, 9)
\item[12.] (0, 1, 4, 7, 4, 9, 5, 8, 1, 2, 6, 9)
\item[13.] (0, 1, 4, 7, 4, 3, 5, 8, 4, 2, 6, 9)
\item[14.] (0, 1, 4, 7, 4, 6, 5, 8, 4, 2, 6, 9)
\item[15.] (0, 1, 4, 7, 4, 9, 5, 8, 4, 2, 6, 9)
\item[16.] (0, 1, 4, 7, 4, 3, 5, 8, 7, 2, 6, 9)
\item[17.] (0, 1, 4, 7, 4, 6, 5, 8, 7, 2, 6, 9)
\item[18.] (0, 1, 4, 7, 4, 9, 5, 8, 7, 2, 6, 9)
\item[19.] (0, 1, 4, 7, 7, 3, 5, 8, 1, 2, 6, 9)
\item[20.] (0, 1, 4, 7, 7, 6, 5, 8, 1, 2, 6, 9)
\item[21.] (0, 1, 4, 7, 7, 9, 5, 8, 1, 2, 6, 9)
\item[22.] (0, 1, 4, 7, 7, 3, 5, 8, 4, 2, 6, 9)
\item[23.] (0, 1, 4, 7, 7, 6, 5, 8, 4, 2, 6, 9)
\item[24.] (0, 1, 4, 7, 7, 9, 5, 8, 4, 2, 6, 9)
\item[25.] (0, 1, 4, 7, 7, 3, 5, 8, 7, 2, 6, 9)
\item[26.] (0, 1, 4, 7, 7, 6, 5, 8, 7, 2, 6, 9)
\item[27.] (0, 1, 4, 7, 7, 9, 5, 8, 7, 2, 6, 9)
\item[28.] (0, 1, 4, 7, 1, 3, 5, 8, 1, 5, 6, 9)
\item[29.] (0, 1, 4, 7, 1, 6, 5, 8, 1, 5, 6, 9)
\item[30.] (0, 1, 4, 7, 1, 9, 5, 8, 1, 5, 6, 9)
\item[31.] (0, 1, 4, 7, 1, 3, 5, 8, 4, 5, 6, 9)
\item[32.] \textbf{(0, 1, 4, 7, 1, 6, 5, 8, 4, 5, 6, 9)}
\item[33.] (0, 1, 4, 7, 1, 9, 5, 8, 4, 5, 6, 9)
\item[34.] (0, 1, 4, 7, 1, 3, 5, 8, 7, 5, 6, 9)
\item[35.] (0, 1, 4, 7, 1, 6, 5, 8, 7, 5, 6, 9)
\item[36.] (0, 1, 4, 7, 1, 9, 5, 8, 7, 5, 6, 9)
\item[37.] (0, 1, 4, 7, 4, 3, 5, 8, 1, 5, 6, 9)
\item[38.] \textbf{(0, 1, 4, 7, 4, 6, 5, 8, 1, 5, 6, 9)}
\item[39.] (0, 1, 4, 7, 4, 9, 5, 8, 1, 5, 6, 9)
\item[40.] (0, 1, 4, 7, 4, 3, 5, 8, 4, 5, 6, 9)
\item[41.] (0, 1, 4, 7, 4, 6, 5, 8, 4, 5, 6, 9)
\item[42.] (0, 1, 4, 7, 4, 9, 5, 8, 4, 5, 6, 9)
\item[43.] (0, 1, 4, 7, 4, 3, 5, 8, 7, 5, 6, 9)
\item[44.] (0, 1, 4, 7, 4, 6, 5, 8, 7, 5, 6, 9)
\item[45.] (0, 1, 4, 7, 4, 9, 5, 8, 7, 5, 6, 9)
\item[46.] (0, 1, 4, 7, 7, 3, 5, 8, 1, 5, 6, 9)
\item[47.] (0, 1, 4, 7, 7, 6, 5, 8, 1, 5, 6, 9)
\item[48.] (0, 1, 4, 7, 7, 9, 5, 8, 1, 5, 6, 9)
\item[49.] (0, 1, 4, 7, 7, 3, 5, 8, 4, 5, 6, 9)
\item[50.] (0, 1, 4, 7, 7, 6, 5, 8, 4, 5, 6, 9)
\item[51.] (0, 1, 4, 7, 7, 9, 5, 8, 4, 5, 6, 9)
\item[52.] (0, 1, 4, 7, 7, 3, 5, 8, 7, 5, 6, 9)
\item[53.] (0, 1, 4, 7, 7, 6, 5, 8, 7, 5, 6, 9)
\item[54.] (0, 1, 4, 7, 7, 9, 5, 8, 7, 5, 6, 9)
\item[55.] (0, 1, 4, 7, 1, 3, 5, 8, 1, 8, 6, 9)
\item[56.] (0, 1, 4, 7, 1, 6, 5, 8, 1, 8, 6, 9)
\item[57.] (0, 1, 4, 7, 1, 9, 5, 8, 1, 8, 6, 9)
\item[58.] (0, 1, 4, 7, 1, 3, 5, 8, 4, 8, 6, 9)
\item[59.] (0, 1, 4, 7, 1, 6, 5, 8, 4, 8, 6, 9)
\item[60.] (0, 1, 4, 7, 1, 9, 5, 8, 4, 8, 6, 9)
\item[61.] (0, 1, 4, 7, 1, 3, 5, 8, 7, 8, 6, 9)
\item[62.] (0, 1, 4, 7, 1, 6, 5, 8, 7, 8, 6, 9)
\item[63.] \textbf{(0, 1, 4, 7, 1, 9, 5, 8, 7, 8, 6, 9)}
\item[64.] (0, 1, 4, 7, 4, 3, 5, 8, 1, 8, 6, 9)
\item[65.] (0, 1, 4, 7, 4, 6, 5, 8, 1, 8, 6, 9)
\item[66.] (0, 1, 4, 7, 4, 9, 5, 8, 1, 8, 6, 9)
\item[67.] (0, 1, 4, 7, 4, 3, 5, 8, 4, 8, 6, 9)
\item[68.] (0, 1, 4, 7, 4, 6, 5, 8, 4, 8, 6, 9)
\item[69.] (0, 1, 4, 7, 4, 9, 5, 8, 4, 8, 6, 9)
\item[70.] (0, 1, 4, 7, 4, 3, 5, 8, 7, 8, 6, 9)
\item[71.] (0, 1, 4, 7, 4, 6, 5, 8, 7, 8, 6, 9)
\item[72.] (0, 1, 4, 7, 4, 9, 5, 8, 7, 8, 6, 9)
\item[73.] (0, 1, 4, 7, 7, 3, 5, 8, 1, 8, 6, 9)
\item[74.] (0, 1, 4, 7, 7, 6, 5, 8, 1, 8, 6, 9)
\item[75.] \textbf{(0, 1, 4, 7, 7, 9, 5, 8, 1, 8, 6, 9)}
\item[76.] (0, 1, 4, 7, 7, 3, 5, 8, 4, 8, 6, 9)
\item[77.] (0, 1, 4, 7, 7, 6, 5, 8, 4, 8, 6, 9)
\item[78.] (0, 1, 4, 7, 7, 9, 5, 8, 4, 8, 6, 9)
\item[79.] (0, 1, 4, 7, 7, 3, 5, 8, 7, 8, 6, 9)
\item[80.] (0, 1, 4, 7, 7, 6, 5, 8, 7, 8, 6, 9)
\item[81.] (0, 1, 4, 7, 7, 9, 5, 8, 7, 8, 6, 9)\hspace{0.3cm}\normalsize$\uparrow\Gamma_1$

\vspace{0.1cm}
      \hrule
      \vspace{-0.1cm}\tiny 
\item[82.] (0, 2, 4, 7, 1, 3, 5, 8, 1, 2, 6, 9)\hspace{0.3cm}\normalsize$\downarrow\Gamma_2$\tiny
\item[83.] (0, 2, 4, 7, 1, 6, 5, 8, 1, 2, 6, 9)
\item[84.] (0, 2, 4, 7, 1, 9, 5, 8, 1, 2, 6, 9)
\item[85.] \textbf{(0, 2, 4, 7, 1, 3, 5, 8, 4, 2, 6, 9)}
\item[86.] (0, 2, 4, 7, 1, 6, 5, 8, 4, 2, 6, 9)
\item[87.] (0, 2, 4, 7, 1, 9, 5, 8, 4, 2, 6, 9)
\item[88.] \textbf{(0, 2, 4, 7, 1, 3, 5, 8, 7, 2, 6, 9)}
\item[89.] (0, 2, 4, 7, 1, 6, 5, 8, 7, 2, 6, 9)
\item[90.] (0, 2, 4, 7, 1, 9, 5, 8, 7, 2, 6, 9)
\item[91.] (0, 2, 4, 7, 4, 3, 5, 8, 1, 2, 6, 9)
\item[92.] (0, 2, 4, 7, 4, 6, 5, 8, 1, 2, 6, 9)
\item[93.] (0, 2, 4, 7, 4, 9, 5, 8, 1, 2, 6, 9)
\item[94.] (0, 2, 4, 7, 4, 3, 5, 8, 4, 2, 6, 9)
\item[95.] (0, 2, 4, 7, 4, 6, 5, 8, 4, 2, 6, 9)
\item[96.] (0, 2, 4, 7, 4, 9, 5, 8, 4, 2, 6, 9)
\item[97.] (0, 2, 4, 7, 4, 3, 5, 8, 7, 2, 6, 9)
\item[98.] (0, 2, 4, 7, 4, 6, 5, 8, 7, 2, 6, 9)
\item[99.] (0, 2, 4, 7, 4, 9, 5, 8, 7, 2, 6, 9)
\item[100.] (0, 2, 4, 7, 7, 3, 5, 8, 1, 2, 6, 9)
\item[101.] (0, 2, 4, 7, 7, 6, 5, 8, 1, 2, 6, 9)
\item[102.] (0, 2, 4, 7, 7, 9, 5, 8, 1, 2, 6, 9)
\item[103.] (0, 2, 4, 7, 7, 3, 5, 8, 4, 2, 6, 9)
\item[104.] (0, 2, 4, 7, 7, 6, 5, 8, 4, 2, 6, 9)
\item[105.] (0, 2, 4, 7, 7, 9, 5, 8, 4, 2, 6, 9)
\item[106.] (0, 2, 4, 7, 7, 3, 5, 8, 7, 2, 6, 9)
\item[107.] (0, 2, 4, 7, 7, 6, 5, 8, 7, 2, 6, 9)
\item[108.] (0, 2, 4, 7, 7, 9, 5, 8, 7, 2, 6, 9)
\item[109.] (0, 2, 4, 7, 1, 3, 5, 8, 1, 5, 6, 9)
\item[110.] (0, 2, 4, 7, 1, 6, 5, 8, 1, 5, 6, 9)
\item[111.] (0, 2, 4, 7, 1, 9, 5, 8, 1, 5, 6, 9)
\item[112.] \textbf{(0, 2, 4, 7, 1, 3, 5, 8, 4, 5, 6, 9)}
\item[113.] \textbf{(0, 2, 4, 7, 1, 6, 5, 8, 4, 5, 6, 9)}
\item[114.] (0, 2, 4, 7, 1, 9, 5, 8, 4, 5, 6, 9)
\item[115.] \textbf{(0, 2, 4, 7, 1, 3, 5, 8, 7, 5, 6, 9)}
\item[116.] (0, 2, 4, 7, 1, 6, 5, 8, 7, 5, 6, 9)
\item[117.] (0, 2, 4, 7, 1, 9, 5, 8, 7, 5, 6, 9)
\item[118.] (0, 2, 4, 7, 4, 3, 5, 8, 1, 5, 6, 9)
\item[119.] \textbf{(0, 2, 4, 7, 4, 6, 5, 8, 1, 5, 6, 9)}
\item[120.] (0, 2, 4, 7, 4, 9, 5, 8, 1, 5, 6, 9)
\item[121.] (0, 2, 4, 7, 4, 3, 5, 8, 4, 5, 6, 9)
\item[122.] (0, 2, 4, 7, 4, 6, 5, 8, 4, 5, 6, 9)
\item[123.] (0, 2, 4, 7, 4, 9, 5, 8, 4, 5, 6, 9)
\item[124.] (0, 2, 4, 7, 4, 3, 5, 8, 7, 5, 6, 9)
\item[125.] (0, 2, 4, 7, 4, 6, 5, 8, 7, 5, 6, 9)
\item[126.] (0, 2, 4, 7, 4, 9, 5, 8, 7, 5, 6, 9)
\item[127.] (0, 2, 4, 7, 7, 3, 5, 8, 1, 5, 6, 9)
\item[128.] (0, 2, 4, 7, 7, 6, 5, 8, 1, 5, 6, 9)
\item[129.] (0, 2, 4, 7, 7, 9, 5, 8, 1, 5, 6, 9)
\item[130.] (0, 2, 4, 7, 7, 3, 5, 8, 4, 5, 6, 9)
\item[131.] (0, 2, 4, 7, 7, 6, 5, 8, 4, 5, 6, 9)
\item[132.] (0, 2, 4, 7, 7, 9, 5, 8, 4, 5, 6, 9)
\item[133.] (0, 2, 4, 7, 7, 3, 5, 8, 7, 5, 6, 9)
\item[134.] (0, 2, 4, 7, 7, 6, 5, 8, 7, 5, 6, 9)
\item[135.] (0, 2, 4, 7, 7, 9, 5, 8, 7, 5, 6, 9)
\item[136.] (0, 2, 4, 7, 1, 3, 5, 8, 1, 8, 6, 9)
\item[137.] (0, 2, 4, 7, 1, 6, 5, 8, 1, 8, 6, 9)
\item[138.] (0, 2, 4, 7, 1, 9, 5, 8, 1, 8, 6, 9)
\item[139.] \textbf{(0, 2, 4, 7, 1, 3, 5, 8, 4, 8, 6, 9)}
\item[140.] (0, 2, 4, 7, 1, 6, 5, 8, 4, 8, 6, 9)
\item[141.] (0, 2, 4, 7, 1, 9, 5, 8, 4, 8, 6, 9)
\item[142.] \textbf{(0, 2, 4, 7, 1, 3, 5, 8, 7, 8, 6, 9)}
\item[143.] (0, 2, 4, 7, 1, 6, 5, 8, 7, 8, 6, 9)
\item[144.] \textbf{(0, 2, 4, 7, 1, 9, 5, 8, 7, 8, 6, 9)}
\item[145.] (0, 2, 4, 7, 4, 3, 5, 8, 1, 8, 6, 9)
\item[146.] (0, 2, 4, 7, 4, 6, 5, 8, 1, 8, 6, 9)
\item[147.] (0, 2, 4, 7, 4, 9, 5, 8, 1, 8, 6, 9)
\item[148.] (0, 2, 4, 7, 4, 3, 5, 8, 4, 8, 6, 9)
\item[149.] (0, 2, 4, 7, 4, 6, 5, 8, 4, 8, 6, 9)
\item[150.] (0, 2, 4, 7, 4, 9, 5, 8, 4, 8, 6, 9)
\item[151.] (0, 2, 4, 7, 4, 3, 5, 8, 7, 8, 6, 9)
\item[152.] (0, 2, 4, 7, 4, 6, 5, 8, 7, 8, 6, 9)
\item[153.] (0, 2, 4, 7, 4, 9, 5, 8, 7, 8, 6, 9)
\item[154.] (0, 2, 4, 7, 7, 3, 5, 8, 1, 8, 6, 9)
\item[155.] (0, 2, 4, 7, 7, 6, 5, 8, 1, 8, 6, 9)
\item[156.] \textbf{(0, 2, 4, 7, 7, 9, 5, 8, 1, 8, 6, 9)}
\item[157.] (0, 2, 4, 7, 7, 3, 5, 8, 4, 8, 6, 9)
\item[158.] (0, 2, 4, 7, 7, 6, 5, 8, 4, 8, 6, 9)
\item[159.] (0, 2, 4, 7, 7, 9, 5, 8, 4, 8, 6, 9)
\item[160.] (0, 2, 4, 7, 7, 3, 5, 8, 7, 8, 6, 9)
\item[161.] (0, 2, 4, 7, 7, 6, 5, 8, 7, 8, 6, 9)
\item[162.] (0, 2, 4, 7, 7, 9, 5, 8, 7, 8, 6, 9)
\item[163.] (0, 5, 4, 7, 1, 3, 5, 8, 1, 2, 6, 9)
\item[164.] (0, 5, 4, 7, 1, 6, 5, 8, 1, 2, 6, 9)
\item[165.] (0, 5, 4, 7, 1, 9, 5, 8, 1, 2, 6, 9)
\item[166.] \textbf{(0, 5, 4, 7, 1, 3, 5, 8, 4, 2, 6, 9)}
\item[167.] \textbf{(0, 5, 4, 7, 1, 6, 5, 8, 4, 2, 6, 9)}
\item[168.] (0, 5, 4, 7, 1, 9, 5, 8, 4, 2, 6, 9)
\item[169.] \textbf{(0, 5, 4, 7, 1, 3, 5, 8, 7, 2, 6, 9)}
\item[170.] (0, 5, 4, 7, 1, 6, 5, 8, 7, 2, 6, 9)
\item[171.] (0, 5, 4, 7, 1, 9, 5, 8, 7, 2, 6, 9)
\item[172.] \textbf{(0, 5, 4, 7, 4, 3, 5, 8, 1, 2, 6, 9)}
\item[173.] \textbf{(0, 5, 4, 7, 4, 6, 5, 8, 1, 2, 6, 9)}
\item[174.] \textbf{(0, 5, 4, 7, 4, 9, 5, 8, 1, 2, 6, 9)}
\item[175.] (0, 5, 4, 7, 4, 3, 5, 8, 4, 2, 6, 9)
\item[176.] (0, 5, 4, 7, 4, 6, 5, 8, 4, 2, 6, 9)
\item[177.] (0, 5, 4, 7, 4, 9, 5, 8, 4, 2, 6, 9)
\item[178.] (0, 5, 4, 7, 4, 3, 5, 8, 7, 2, 6, 9)
\item[179.] (0, 5, 4, 7, 4, 6, 5, 8, 7, 2, 6, 9)
\item[180.] (0, 5, 4, 7, 4, 9, 5, 8, 7, 2, 6, 9)
\item[181.] \textbf{(0, 5, 4, 7, 7, 3, 5, 8, 1, 2, 6, 9)}
\item[182.] \textbf{(0, 5, 4, 7, 7, 6, 5, 8, 1, 2, 6, 9)}
\item[183.] \textbf{(0, 5, 4, 7, 7, 9, 5, 8, 1, 2, 6, 9)}
\item[184.] (0, 5, 4, 7, 7, 3, 5, 8, 4, 2, 6, 9)
\item[185.] (0, 5, 4, 7, 7, 6, 5, 8, 4, 2, 6, 9)
\item[186.] (0, 5, 4, 7, 7, 9, 5, 8, 4, 2, 6, 9)
\item[187.] (0, 5, 4, 7, 7, 3, 5, 8, 7, 2, 6, 9)
\item[188.] (0, 5, 4, 7, 7, 6, 5, 8, 7, 2, 6, 9)
\item[189.] (0, 5, 4, 7, 7, 9, 5, 8, 7, 2, 6, 9)
\item[190.] (0, 5, 4, 7, 1, 3, 5, 8, 1, 5, 6, 9)
\item[191.] (0, 5, 4, 7, 1, 6, 5, 8, 1, 5, 6, 9)
\item[192.] (0, 5, 4, 7, 1, 9, 5, 8, 1, 5, 6, 9)
\item[193.] \textbf{(0, 5, 4, 7, 1, 3, 5, 8, 4, 5, 6, 9)}
\item[194.] (0, 5, 4, 7, 1, 6, 5, 8, 4, 5, 6, 9)
\item[195.] (0, 5, 4, 7, 1, 9, 5, 8, 4, 5, 6, 9)
\item[196.] \textbf{(0, 5, 4, 7, 1, 3, 5, 8, 7, 5, 6, 9)}
\item[197.] (0, 5, 4, 7, 1, 6, 5, 8, 7, 5, 6, 9)
\item[198.] (0, 5, 4, 7, 1, 9, 5, 8, 7, 5, 6, 9)
\item[199.] (0, 5, 4, 7, 4, 3, 5, 8, 1, 5, 6, 9)
\item[200.] (0, 5, 4, 7, 4, 6, 5, 8, 1, 5, 6, 9)
\item[201.] (0, 5, 4, 7, 4, 9, 5, 8, 1, 5, 6, 9)
\item[202.] \textbf{(0, 5, 4, 7, 4, 3, 5, 8, 4, 5, 6, 9)}
\item[203.] \textbf{(0, 5, 4, 7, 4, 6, 5, 8, 4, 5, 6, 9)}
\item[204.] \textbf{(0, 5, 4, 7, 4, 9, 5, 8, 4, 5, 6, 9)}
\item[205.] (0, 5, 4, 7, 4, 3, 5, 8, 7, 5, 6, 9)
\item[206.] (0, 5, 4, 7, 4, 6, 5, 8, 7, 5, 6, 9)
\item[207.] (0, 5, 4, 7, 4, 9, 5, 8, 7, 5, 6, 9)
\item[208.] (0, 5, 4, 7, 7, 3, 5, 8, 1, 5, 6, 9)
\item[209.] (0, 5, 4, 7, 7, 6, 5, 8, 1, 5, 6, 9)
\item[210.] (0, 5, 4, 7, 7, 9, 5, 8, 1, 5, 6, 9)
\item[211.] (0, 5, 4, 7, 7, 3, 5, 8, 4, 5, 6, 9)
\item[212.] (0, 5, 4, 7, 7, 6, 5, 8, 4, 5, 6, 9)
\item[213.] (0, 5, 4, 7, 7, 9, 5, 8, 4, 5, 6, 9)
\item[214.] (0, 5, 4, 7, 7, 3, 5, 8, 7, 5, 6, 9)
\item[215.] (0, 5, 4, 7, 7, 6, 5, 8, 7, 5, 6, 9)
\item[216.] (0, 5, 4, 7, 7, 9, 5, 8, 7, 5, 6, 9)
\item[217.] (0, 5, 4, 7, 1, 3, 5, 8, 1, 8, 6, 9)
\item[218.] (0, 5, 4, 7, 1, 6, 5, 8, 1, 8, 6, 9)
\item[219.] (0, 5, 4, 7, 1, 9, 5, 8, 1, 8, 6, 9)
\item[220.] \textbf{(0, 5, 4, 7, 1, 3, 5, 8, 4, 8, 6, 9)}
\item[221.] (0, 5, 4, 7, 1, 6, 5, 8, 4, 8, 6, 9)
\item[222.] (0, 5, 4, 7, 1, 9, 5, 8, 4, 8, 6, 9)
\item[223.] \textbf{(0, 5, 4, 7, 1, 3, 5, 8, 7, 8, 6, 9)}
\item[224.] (0, 5, 4, 7, 1, 6, 5, 8, 7, 8, 6, 9)
\item[225.] (0, 5, 4, 7, 1, 9, 5, 8, 7, 8, 6, 9)
\item[226.] (0, 5, 4, 7, 4, 3, 5, 8, 1, 8, 6, 9)
\item[227.] (0, 5, 4, 7, 4, 6, 5, 8, 1, 8, 6, 9)
\item[228.] (0, 5, 4, 7, 4, 9, 5, 8, 1, 8, 6, 9)
\item[229.] (0, 5, 4, 7, 4, 3, 5, 8, 4, 8, 6, 9)
\item[230.] (0, 5, 4, 7, 4, 6, 5, 8, 4, 8, 6, 9)
\item[231.] (0, 5, 4, 7, 4, 9, 5, 8, 4, 8, 6, 9)
\item[232.] (0, 5, 4, 7, 4, 3, 5, 8, 7, 8, 6, 9)
\item[233.] \textbf{(0, 5, 4, 7, 4, 6, 5, 8, 7, 8, 6, 9)}
\item[234.] \textbf{(0, 5, 4, 7, 4, 9, 5, 8, 7, 8, 6, 9)}
\item[235.] (0, 5, 4, 7, 7, 3, 5, 8, 1, 8, 6, 9)
\item[236.] (0, 5, 4, 7, 7, 6, 5, 8, 1, 8, 6, 9)
\item[237.] (0, 5, 4, 7, 7, 9, 5, 8, 1, 8, 6, 9)
\item[238.] (0, 5, 4, 7, 7, 3, 5, 8, 4, 8, 6, 9)
\item[239.] \textbf{(0, 5, 4, 7, 7, 6, 5, 8, 4, 8, 6, 9)}
\item[240.] \textbf{(0, 5, 4, 7, 7, 9, 5, 8, 4, 8, 6, 9)}
\item[241.] (0, 5, 4, 7, 7, 3, 5, 8, 7, 8, 6, 9)
\item[242.] (0, 5, 4, 7, 7, 6, 5, 8, 7, 8, 6, 9)
\item[243.] (0, 5, 4, 7, 7, 9, 5, 8, 7, 8, 6, 9)
\item[244.] (0, 8, 4, 7, 1, 3, 5, 8, 1, 2, 6, 9)
\item[245.] (0, 8, 4, 7, 1, 6, 5, 8, 1, 2, 6, 9)
\item[246.] (0, 8, 4, 7, 1, 9, 5, 8, 1, 2, 6, 9)
\item[247.] \textbf{(0, 8, 4, 7, 1, 3, 5, 8, 4, 2, 6, 9)}
\item[248.] (0, 8, 4, 7, 1, 6, 5, 8, 4, 2, 6, 9)
\item[249.] (0, 8, 4, 7, 1, 9, 5, 8, 4, 2, 6, 9)
\item[250.] \textbf{(0, 8, 4, 7, 1, 3, 5, 8, 7, 2, 6, 9)}
\item[251.] (0, 8, 4, 7, 1, 6, 5, 8, 7, 2, 6, 9)
\item[252.] \textbf{(0, 8, 4, 7, 1, 9, 5, 8, 7, 2, 6, 9)}
\item[253.] \textbf{(0, 8, 4, 7, 4, 3, 5, 8, 1, 2, 6, 9)}
\item[254.] \textbf{(0, 8, 4, 7, 4, 6, 5, 8, 1, 2, 6, 9)}
\item[255.] \textbf{(0, 8, 4, 7, 4, 9, 5, 8, 1, 2, 6, 9)}
\item[256.] (0, 8, 4, 7, 4, 3, 5, 8, 4, 2, 6, 9)
\item[257.] (0, 8, 4, 7, 4, 6, 5, 8, 4, 2, 6, 9)
\item[258.] (0, 8, 4, 7, 4, 9, 5, 8, 4, 2, 6, 9)
\item[259.] (0, 8, 4, 7, 4, 3, 5, 8, 7, 2, 6, 9)
\item[260.] (0, 8, 4, 7, 4, 6, 5, 8, 7, 2, 6, 9)
\item[261.] (0, 8, 4, 7, 4, 9, 5, 8, 7, 2, 6, 9)
\item[262.] \textbf{(0, 8, 4, 7, 7, 3, 5, 8, 1, 2, 6, 9)}
\item[263.] \textbf{(0, 8, 4, 7, 7, 6, 5, 8, 1, 2, 6, 9)}
\item[264.] \textbf{(0, 8, 4, 7, 7, 9, 5, 8, 1, 2, 6, 9)}
\item[265.] (0, 8, 4, 7, 7, 3, 5, 8, 4, 2, 6, 9)
\item[266.] (0, 8, 4, 7, 7, 6, 5, 8, 4, 2, 6, 9)
\item[267.] (0, 8, 4, 7, 7, 9, 5, 8, 4, 2, 6, 9)
\item[268.] (0, 8, 4, 7, 7, 3, 5, 8, 7, 2, 6, 9)
\item[269.] (0, 8, 4, 7, 7, 6, 5, 8, 7, 2, 6, 9)
\item[270.] (0, 8, 4, 7, 7, 9, 5, 8, 7, 2, 6, 9)
\item[271.] (0, 8, 4, 7, 1, 3, 5, 8, 1, 5, 6, 9)
\item[272.] (0, 8, 4, 7, 1, 6, 5, 8, 1, 5, 6, 9)
\item[273.] (0, 8, 4, 7, 1, 9, 5, 8, 1, 5, 6, 9)
\item[274.] \textbf{(0, 8, 4, 7, 1, 3, 5, 8, 4, 5, 6, 9)}
\item[275.] (0, 8, 4, 7, 1, 6, 5, 8, 4, 5, 6, 9)
\item[276.] (0, 8, 4, 7, 1, 9, 5, 8, 4, 5, 6, 9)
\item[277.] \textbf{(0, 8, 4, 7, 1, 3, 5, 8, 7, 5, 6, 9)}
\item[278.] (0, 8, 4, 7, 1, 6, 5, 8, 7, 5, 6, 9)
\item[279.] (0, 8, 4, 7, 1, 9, 5, 8, 7, 5, 6, 9)
\item[280.] (0, 8, 4, 7, 4, 3, 5, 8, 1, 5, 6, 9)
\item[281.] (0, 8, 4, 7, 4, 6, 5, 8, 1, 5, 6, 9)
\item[282.] (0, 8, 4, 7, 4, 9, 5, 8, 1, 5, 6, 9)
\item[283.] (0, 8, 4, 7, 4, 3, 5, 8, 4, 5, 6, 9)
\item[284.] (0, 8, 4, 7, 4, 6, 5, 8, 4, 5, 6, 9)
\item[285.] (0, 8, 4, 7, 4, 9, 5, 8, 4, 5, 6, 9)
\item[286.] (0, 8, 4, 7, 4, 3, 5, 8, 7, 5, 6, 9)
\item[287.] \textbf{(0, 8, 4, 7, 4, 6, 5, 8, 7, 5, 6, 9)}
\item[288.] \textbf{(0, 8, 4, 7, 4, 9, 5, 8, 7, 5, 6, 9)}
\item[289.] (0, 8, 4, 7, 7, 3, 5, 8, 1, 5, 6, 9)
\item[290.] (0, 8, 4, 7, 7, 6, 5, 8, 1, 5, 6, 9)
\item[291.] (0, 8, 4, 7, 7, 9, 5, 8, 1, 5, 6, 9)
\item[292.] (0, 8, 4, 7, 7, 3, 5, 8, 4, 5, 6, 9)
\item[293.] \textbf{(0, 8, 4, 7, 7, 6, 5, 8, 4, 5, 6, 9)}
\item[294.] \textbf{(0, 8, 4, 7, 7, 9, 5, 8, 4, 5, 6, 9)}
\item[295.] (0, 8, 4, 7, 7, 3, 5, 8, 7, 5, 6, 9)
\item[296.] (0, 8, 4, 7, 7, 6, 5, 8, 7, 5, 6, 9)
\item[297.] (0, 8, 4, 7, 7, 9, 5, 8, 7, 5, 6, 9)
\item[298.] (0, 8, 4, 7, 1, 3, 5, 8, 1, 8, 6, 9)
\item[299.] (0, 8, 4, 7, 1, 6, 5, 8, 1, 8, 6, 9)
\item[300.] (0, 8, 4, 7, 1, 9, 5, 8, 1, 8, 6, 9)
\item[301.] \textbf{(0, 8, 4, 7, 1, 3, 5, 8, 4, 8, 6, 9)}
\item[302.] (0, 8, 4, 7, 1, 6, 5, 8, 4, 8, 6, 9)
\item[303.] (0, 8, 4, 7, 1, 9, 5, 8, 4, 8, 6, 9)
\item[304.] \textbf{(0, 8, 4, 7, 1, 3, 5, 8, 7, 8, 6, 9)}
\item[305.] (0, 8, 4, 7, 1, 6, 5, 8, 7, 8, 6, 9)
\item[306.] (0, 8, 4, 7, 1, 9, 5, 8, 7, 8, 6, 9)
\item[307.] (0, 8, 4, 7, 4, 3, 5, 8, 1, 8, 6, 9)
\item[308.] (0, 8, 4, 7, 4, 6, 5, 8, 1, 8, 6, 9)
\item[309.] (0, 8, 4, 7, 4, 9, 5, 8, 1, 8, 6, 9)
\item[310.] (0, 8, 4, 7, 4, 3, 5, 8, 4, 8, 6, 9)
\item[311.] (0, 8, 4, 7, 4, 6, 5, 8, 4, 8, 6, 9)
\item[312.] (0, 8, 4, 7, 4, 9, 5, 8, 4, 8, 6, 9)
\item[313.] (0, 8, 4, 7, 4, 3, 5, 8, 7, 8, 6, 9)
\item[314.] (0, 8, 4, 7, 4, 6, 5, 8, 7, 8, 6, 9)
\item[315.] (0, 8, 4, 7, 4, 9, 5, 8, 7, 8, 6, 9)
\item[316.] (0, 8, 4, 7, 7, 3, 5, 8, 1, 8, 6, 9)
\item[317.] (0, 8, 4, 7, 7, 6, 5, 8, 1, 8, 6, 9)
\item[318.] (0, 8, 4, 7, 7, 9, 5, 8, 1, 8, 6, 9)
\item[319.] (0, 8, 4, 7, 7, 3, 5, 8, 4, 8, 6, 9)
\item[320.] (0, 8, 4, 7, 7, 6, 5, 8, 4, 8, 6, 9)
\item[321.] (0, 8, 4, 7, 7, 9, 5, 8, 4, 8, 6, 9)
\item[322.] \textbf{(0, 8, 4, 7, 7, 3, 5, 8, 7, 8, 6, 9)}
\item[323.] \textbf{(0, 8, 4, 7, 7, 6, 5, 8, 7, 8, 6, 9)}
\item[324.] \textbf{(0, 8, 4, 7, 7, 9, 5, 8, 7, 8, 6, 9)}
\end{itemize}
\end{multicols}

\normalsize 

The following are the index numbers of the 123 encodings in the 4D expansion of the 2D sunflower which provide Formality graphs $\Gamma$ whose formulas $\phi(\Gamma)$ are linearly independent. Note that indices start from 1. In bold are the indices of the 64 4D sunflower encodings which provide Formality graphs $\Gamma$ whose skew pairs $\frac{1}{2}\Bigl(\phi\bigl(\Gamma(a^1,a^2)\bigr)-\phi\bigl(\Gamma(a^2,a^1)\bigr)\Bigr)$ are linearly independent. 

\bigskip 

\textbf{1}, \textbf{2}, 3, \textbf{4}, \textbf{5}, \textbf{6}, 7, 8, 9, \textbf{11}, \textbf{12}, \textbf{13}, \textbf{14}, \textbf{15}, \textbf{16}, \textbf{17}, \textbf{18}, 20, 21, 23, 24, 25, 26, 27, \textbf{29}, \textbf{30}, \textbf{33}, \textbf{35}, \textbf{36}, 39, \textbf{41}, \textbf{42}, \textbf{44}, \textbf{45}, 48, \textbf{51}, \textbf{53}, 54, 57, 60, 69, 72, 81, \textbf{82}, \textbf{83}, 84, \textbf{94}, \textbf{95}, \textbf{96}, \textbf{97}, \textbf{98}, 103, 105, 106, 107, 108, \textbf{109}, \textbf{110}, \textbf{111}, \textbf{114}, \textbf{117}, \textbf{118}, 120, \textbf{127}, 129, 136, 138, 145, 154, \textbf{163}, \textbf{165}, \textbf{168}, 171, \textbf{175}, \textbf{177}, 180, \textbf{186}, \textbf{187}, \textbf{188}, \textbf{189}, \textbf{190}, \textbf{191}, \textbf{192}, 195, \textbf{197}, \textbf{198}, \textbf{206}, \textbf{207}, \textbf{208}, 210, \textbf{211}, \textbf{213}, \textbf{214}, \textbf{215}, \textbf{216}, \textbf{217}, \textbf{218}, 219, 224, 236, \textbf{241}, \textbf{242}, 243, 244, 245, 256, 268, 269, 271, 283, 285, 296, 298, 299, 300, 302, 303, 307, 310, 311, 312, 313, 314.

\bigskip 

For clarity, below are the indices of the 64 encodings in the 4D expansion of the 2D sunflower which provide Formality graphs $\Gamma$ whose skew pairs~$\frac{1}{2}\Bigl(\phi\bigl(\Gamma(a^1,a^2)\bigr)-\phi\bigl(\Gamma(a^2,a^1)\bigr)\Bigr)$ are linearly independent. Note that indices start from 1.

\bigskip 

1, 2, 4, 5, 6, 11, 12, 13, 14, 15, 16, 17, 18, 29, 30, 33, 35, 36, 41, 42, 44, 45, 51, 53, 82, 83, 94, 95, 96, 97, 98, 109, 110, 111, 114, 117, 118, 127, 163, 165, 168, 175, 177, 186, 187, 188, 189, 190, 191, 192, 197, 198, 206, 207, 208, 211, 213, 214, 215, 216, 217, 218, 241, 242.

\end{document}